\documentclass[a4paper,10pt]{amsart}
\usepackage{amsmath,amsfonts,amssymb,amsthm, xypic,mathrsfs}
\usepackage{epsfig,psfrag}
\usepackage[all]{xy}

\def\Z{\mathbb{Z}}

\def\CC{\mathbb{C}}

\def\H{\mathfrak{h}}

\DeclareMathOperator{\Coh}{Coh}
\DeclareMathOperator{\Rep}{Rep}

\newtheorem{theo}{\bf{Theorem}}[section]
\newtheorem{lem}[theo]{Lemma}

\newtheorem{prop}[theo]{Proposition}

\newtheorem{conj}[theo]{Conjecture}

\title{A new involution for quantum loop algebras}

\author{Jyun-Ao Lin}

\begin{document}

\begin{abstract}
In this article, we introduce a completion $\widehat{U}^+_v(\mathcal{L}\mathfrak{g})$ of the positive half of the quantum affinization $U^+_v(\mathcal{L}\mathfrak{g})$ of a symmetrizable Kac-Moody algebra $\mathfrak{g}$. On $\widehat{U}^+_v(\mathcal{L}(\mathfrak{g}))$, we define a new \lq\lq bar-involution \rq\rq \; and construct the analogue Kashiwara's operators. We conjecture that the resulting pair $(\widehat{\mathcal{L}},\widehat{\mathcal{B}})$ is a crystal basis which provides the existence of the \lq\lq canonical basis \rq\rq \; on the (completion of the) of the positive half of the quamtum affinization. 
\end{abstract}
\maketitle
\setcounter{section}{-1}

\section{Introduction}
One of the main outcomes of the theory of quantum groups is the discovery, by Kashiwara \cite{Kash91} and Lusztig \cite{Lu90}, of the canonical bases in quantized enveloping algebras with certain favorable properties: positivity of structure constants, compatibility with all highest weight integrable representations, etc. These canonical bases have been proven to be powerful tool in the study of the representation theory of quantum Kac-Moody algebras, encoding character formulas and decomposition numbers.

\vspace{.1in}

The construction of the canonical basis $\mathbf{B}$ of a quantized enveloping algebra $U^+_v(\mathfrak{g})$ in \cite{Lu90} is based on Ringel's discovery in 90s that $U^+_v(\mathfrak{g})$ can be realized in the Hall algebra of the category of representations of a quiver $Q$ with underlying graph is the Dynkin diagram of $\mathfrak{g}$ (\cite{Ring90}). The set of isomorphism classes of representations of $Q$ of a given class $\mathbf{d} \in K_0(\Rep Q)$ is the set of orbits of a reductive group $G_{\mathbf{d}}$ on a vector sapce $E_{\mathbf{d}}$ and Lusztig realizes $U_v^+(\mathfrak{g})$ geometrically as a convolution algebra of semisimple, $G_{\mathbf{d}}$-equivariant constructible sheaves on $E_{\mathbf{d}}$ and obtains the canonical basis as the set of all simple perverse sheaves on this algebra. In \cite{Kap97}, Kapranov shows that the Hall algebra of the category of coherent sheaves on a smooth projective line provides a realization of the affinization(c.f. \cite{Her05}) $U^+_v(\mathcal{L}(\mathfrak{sl}_2))$ of the Drinfeld's positive part of the quantum affine algebra $U_v(\widehat{\mathfrak{sl}}_2)$. Schiffmann then constructs in \cite{Sch06} the canonical basis $\widehat{\mathcal{B}}$ of the completion $\widehat{U}^+_v(\mathcal{L}\mathfrak{sl}_2)$ in a similar fashion and prove the compatibility with some integrable lowest weight representations. 

\vspace{.1in}

Kashiwara's scheme to construct the canonical basis is quite different from Lusztig's one and his approach makes sense for all symmetriable Kac-Moody algebra $\mathfrak{g}$. The main ingredients of \cite{Kash91} are certain operators(called Kashiwara's operators) $\widetilde{E}_i,\widetilde{F}_i : U^+_v(\mathfrak{g}) \to U^+_v(\mathfrak{g})$ for all $i \in I$ to generate the $\mathcal{A}$-lattice $\mathcal{L}$ of $U^+_v(\mathfrak{g})$ and the basis $\mathcal{B}$ of $\mathcal{L} / v \mathcal{L}$, where $\mathcal{A}$ is the localization of $\mathbb{Q}[v]$ at $v=0$. Such a pair $(\mathcal{L},\mathcal{B})$ stable by Kashiwara's operators is called a crystal basis of $U^+_v(\mathfrak{g})$. Consider the so-called bar-involution $\varphi : U^+_v(\mathfrak{g}) \to U^+_v(\mathfrak{g})$ defined by $v \mapsto v^{-1}, E_i \mapsto E_i$ and let $\mathcal{L}^- = \varphi(\mathcal{L})$. Then there is an isomorphism $U^+_v(\mathfrak{g}) \cap \mathcal{L} \cap \mathcal{L}^- \simeq \mathcal{L} / v \mathcal{L}$ and the pre-image of $\mathcal{B}$ under the isomorphism is a basis of $U^+_v(\mathfrak{g})$ which coincides with the canonical basis obtained by Lusztig. Such a triple $(U^+_v(\mathfrak{g}), \mathcal{L}, \mathcal{L}^-)$ equipped with the above isomorphism is called a \textit{balanced triple} and its existence is equivalent to the existence of the canonical basis. 

\vspace{.1in}

In this short paper, we develop a purely algebraic approach (under the scheme of Kashiwara in \cite{Kash91}) to Schiffmann's canonical basis $\widehat{\mathcal{B}}$ on  $\widehat{U}^+_v(\mathcal{L}\mathfrak{sl}_2)$ and extend the construction to the (positive part of the) quantum affinization of any Kac-Moody algebra: We provide a construction of Kashiwara's operators adapted to the Drinfeld's half part of the quantum affinization $U_v(\mathcal{L}\mathfrak{g})$ for all symmetrizable Kac-Moody algebra $\mathfrak{g}$ and use them to generate the (conjectural) crystal basis $(\mathcal{L},\mathcal{B})$. To generalize the concept of the canonical basis, we extend the bar-involution $\varphi$ induced by Verdier duality in the context of $\widehat{U}^+_v(\mathcal{L}\mathfrak{sl}_2)$ to the general cases which seems to be unknown before in the rich study of quantum affine algebras. Since the image of the involution $\varphi$ involves infinite sums, we introduce a certain completion $\widehat{U}^+_v(\mathcal{L}\mathfrak{g})$ of $U^+_v(\mathcal{L}\mathfrak{g})$ and the resulting (conjectural) canonical basis $\widehat{\mathcal{B}}$(if it exists) should actually lie in the completion. 

\vspace{.1in}

Unfortunately the construction of Kashiwara's operators can not apply to the (highest $l$-weight) representations since there is not known non-degenerate bilinear form on them and therefore the author fails to generalize the \lq\lq Grand Loop \rq\rq \; argument to the quantum affinizations. In \cite{Joseph}, Joseph's refinement of Grand Loop argument relies on the representation theory of quantum Weyl algebras which, in our cases, there is an \lq\lq affinization \rq\rq \; of quantum Weyl algebra whose representations seems worth to study. On the other hand, the Grand Loop argument relies heavily on the tensor product of integrable highest weight representations. Here, we might instead by consider the \lq\lq fusion product \rq\rq \; of highest $l$-weight modules(cf. \cite{Her05}). We hope to be able to say more about these problems in the future.

\section{Quantum loop algebras}

\subsection{Quantum groups}
Let $q$ be indeterminate and set $v = q^{-1/2}$. For nonnegative integers $l \geq 0$, define
\begin{equation*}
[l]_{v} := \frac{v^{-l}-v^{l}}{v^{-1}-v}, [l]_{v}! :=[l]_{v}[l-1]_{v} \cdots [2]_{v}[1]_{v}
\end{equation*}

\vspace{.1in}

Let $A=(a_{i,j})_{1 
\leq i,j \leq n}$ be a symmetrizable generalized Cartan matrix, that is, $a_{i,j} \in \Z$, $a_{i,i} = 2$, $a_{i,j} \leq 0$ if $i \neq j$, $a_{i,j} = 0 \iff a_{j,i} = 0$ and there is a matrix $D = \mathrm{diag}(r_1, \dots, r_n)$ with $r_i \in \Z_{>0}$ such that $B = DA = (b_{i,j})_{1 \leq i,j \leq n}$ is symmetric. We denote $I = \{1, \dots, n \}$ and $n = \mathrm{rank}(A)$. We consider a realization $( \H, \H^*, \Pi, \Pi^\vee )$ of $A$: $\H$ is a vector space of dimension $2n - \mathrm{rank}(A)$, $\H^*$ its dual , set of simple roots $\Pi = \{\alpha_1, \dots, \alpha_n \} \subset \H^*$ and set of simple coroots $\Pi^\vee = \{\alpha^\vee_1, \dots, \alpha^\vee_n \} \subset \H$ such that $\alpha_j(\alpha^\vee_i) = a_{i,j}$ for $i,j \in I$. Denote by $\omega_1, \dots, \omega_n \in \H^*$ the fundamental weights. Let $( \ , \ )$ be a nondegenerate bilinear form on $\H^*$ satisfying $( \alpha_i, h) = h(r_i \alpha^\vee_i)$ and $\nu : \H^* \to \H$ be the induced isomorphism. We have in particular $\nu (\alpha_i) = r_i \alpha^\vee_i$ and for any $\lambda, \mu \in \H^*$, $\lambda(\nu(\mu)) = \mu(\nu(\lambda))$.

\vspace{.1in}

We denote by $P = \{ \lambda \in \H^* \mid \forall i \in I, \lambda(\alpha^\vee_i) \in \Z  \}$ the set of weights and $P^+ = \{\lambda \in P | \forall i \in I, \lambda(\alpha^\vee_i) \geq 0 \}$ the set of dominant weights. Let $Q = \bigoplus_{i \in I} \Z \alpha_i \subset P$ be the root lattice and $Q^+ = \sum_{i \in I} \Z_{\geq 0}\alpha_i \subset Q$. For $\lambda, \mu \in \H^*$, we define $\lambda \geq \mu$ if $\lambda - \mu \in Q^+$. Let $\widehat{Q} = Q \times \Z$ and similarly $\widehat{Q}^+ = Q^+ \times \Z$. For $\alpha = \sum_{i \in I} n_i \alpha_i \in Q^+$, we set $|\alpha| = \sum_{i \in I} |n_i|$ and set $Q^+(l) = \{ \alpha \in Q^+ | |\alpha| \leq l \}$. Finally set $\widehat{Q}^+(l) = Q^+(l) \times \Z$.

\vspace{.1in}

Let $\mathfrak{g}$ be the complex Kac-Moody algebra associated to $A$. The quantum group $U_v(\mathfrak{g})$ is the Hopf algebra over $\CC(v)$ generated by elements $E_i, F_i$ for $i \in I$ and $K_h$ for $h \in \H$ subject to the relations

$$K_h K_{h'} = K_{h+h'}, \; K_0 = 1,$$
$$K_h E_j K_{-h} = v^{-\alpha_j(h)}E_j $$
$$K_h F_j K_{-h} = v^{\alpha_j(h)}F_j$$
$$[E_i, F_j ] = \delta_{i,j}\frac{K_{r_i \alpha^\vee_i} - K_{-r_i \alpha^\vee_i}}{v_i^{-1} - v_i}$$
\begin{equation*}
\begin{split}
&\sum_{k=0}^{r} (-1)^k \begin{bmatrix}
r \\ k
\end{bmatrix}_{v_i} E_i^{r-k}E_j E_i^k = 0 \\
&\sum_{k=0}^{r} (-1)^k \begin{bmatrix}
r \\ k
\end{bmatrix}_{v_i} F_i^{r-k}F_j F_i^k = 0
\end{split}
\end{equation*}
for all $i \neq j$ and $r = 1 - a_{i,j}$, where $v_i = v^{r_i}$
The coproduct is given by the formulas
\begin{equation*}
\Delta(K_h) = K_h \otimes K_h, \; \Delta(E_i) = E_i \otimes 1 + K_i \otimes E_i, \; \Delta(F_i) = 1 \otimes F_i + F_i \otimes K_i^{-1},
\end{equation*}
where $K_i = K_{r_i \alpha_i^\vee}$, and the antipode is 
\begin{equation*}
S(K_h) = K_{-h}, \; S(E_i) = - K_i E_i, \; S(F_i) = - F_iK_i^{-1}.
\end{equation*}

\subsection{Quantum loop algebras}
The quantum loop algebra $U_v(\mathcal{L}\mathfrak{g})$ is the $\CC(v)$-algebra generated by $E_{i,l}, F_{i,l}$, for $i \in I$, $l \in \Z$, $K_h$ for $h \in \H$ and $H_{i,s}$ for $1 \leq i \leq n$, $s \in \Z^*$ subject to the following set of relations:

$$K_h K_{h'} = K_{h+h'}, \; K_0 = 1,$$
$$[K_h, H_{j,s}] = 0, \; [H_{i,s},H_{j,t}] = 0, $$
$$K_h E_{j,l} K_{-h} = v^{-\alpha_j(h)}E_{j,l}, \qquad K_h F_{j,l} K_{-h} = v^{\alpha_j(h)}F_{j,l}, $$ 
$$[H_{i,s},E_{j,l}] = \frac{[sb_{i,j}]_v}{s}E_{j,s+l}, $$
$$[H_{i,s},F_{j,l}] = -\frac{[sb_{i,j}]_v}{s}F_{j,s+l},$$
$$[E_{i,l},F_{j,m}] = \delta_{i,j} \frac{\phi^+_{i,l+m} - \phi^-_{i,l+m}}{v_i^{-1} - v_i},$$
\begin{equation}\label{relation1}
\begin{split}
v^{b_{i,j}}E_{i,l+1}E_{j,m} - E_{j,m}E_{i,l+1} &= E_{i,l}E_{j,m+1} - v^{b_{i,j}}E_{j,m+1}E_{i,l}, \\
v^{-b_{i,j}}F_{i,l+1}F_{j,m} - F_{j,m}F_{i,l+1} &= F_{i,l}F_{j,m+1} - v^{-b_{i,j}}F_{j,m+1}F_{i,l},
\end{split}
\end{equation}
\begin{equation}\label{serre}
\begin{split}
\mathrm{For\;}i \neq j,& \; r = 1-a_{i,j}  \; \mathrm{and\; all\; sequences\; of\; integers\; } l_1 \cdots l_r,\\
&\sum_{\sigma \in \mathfrak{S}_r} \sum_{k=0}^r (-1)^k \begin{bmatrix}
r \\ k
\end{bmatrix}_{v_i} 
E_{i,l_{\sigma(1)}} \cdots E_{i,l_{\sigma(k)}}E_{j,l'}E_{i,l_{\sigma(k+1)}} \cdots E_{i,l_{\sigma(r)}} = 0 \\
&\sum_{\sigma \in \mathfrak{S}_r} \sum_{k=0}^r (-1)^k \begin{bmatrix}
r \\ k
\end{bmatrix}_{v_i} 
F_{i,l_{\sigma(1)}} \cdots F_{i,l_{\sigma(k)}}F_{j,l'}F_{i,l_{\sigma(k+1)}} \cdots F_{i,l_{\sigma(r)}} = 0,
\end{split}
\end{equation}
where $\mathfrak{S}_r$ is the symmetric group on $r$ letters and where $\phi^{\pm}_{i,s}$'s are defined by the following equations:
\begin{equation*}
\sum_{s \geq 0} \phi^{\pm}_{i,\pm s} z^{\pm s} = K_{\pm r_i \alpha_i^\vee} \mathrm{exp}\big( \pm (v^{-1} - v) \sum_{s=1}^\infty H_{i, \pm s} z^{\pm s} \big),
\end{equation*}
and $\phi^{\pm}_{i,\pm s} = 0$ for $s < 0$. The relations (\ref{serre}) are the loop analogs of the Serre relations in quantum groups. 

\subsection{Relations between the currents}
It is also convenient to write the defining relations of $U_v(\mathcal{L}\mathfrak{g})$ in terms of formal generating functions(also called currents):
$$E_i(z) = \sum_{l \in \Z} E_{i,l} z^l, \; F_i(z) = \sum_{l \in \Z} F_{i,l} z^l, \; \phi^{\pm}_i(z) = \sum_{s \geq 0} \phi^{\pm}_{i,\pm s} z^{\pm s}$$
and the relations have the form:
$$K_h \phi^{\pm}_j(z) = \phi^{\pm}_j(z)K_h,$$
$$K_h E_j(z) = v^{-\alpha_j(h)} E_j(z)K_h,$$
$$K_h F_j(z) = v^{\alpha_j(h)} F_j(z)K_h,$$
\begin{equation}\label{loop1}
\phi^{\pm}_i(z)E_j(w) = \frac{1-\frac{z}{w}v^{ b_{i,j}}}{v^{ b_{i,j}} - \frac{z}{w}}E_j(w)\phi^{\pm}_i(z),
\end{equation}
$$
\phi^{\pm}_i(z)F_j(w) = \frac{1-\frac{z}{w}v^{ -b_{i,j}}}{v^{ -b_{i,j}} - \frac{z}{w}}F_j(w)\phi^{\pm}_i(z),$$
$$[E_i(z),F_j(w)] = \frac{\delta_{i,j}}{v_i^{-1}-v_i}\{ \delta(\frac{w}{z})\phi^+_i(w) - \delta(\frac{z}{w})\phi^-_i(z),$$
\begin{equation}\label{loop2}
(wv^{ b_{i,j}} - z)E_i(z)E_j(w) = (w-zv^{ b_{i,j}})E_j(w)E_i(z),
\end{equation}
$$(wv^{ -b_{i,j}} - z)F_i(z)F_j(w) = (w-zv^{ -b_{i,j}})F_j(w)F_i(z),$$
\begin{equation}\label{loopserre1}
\sum_{\sigma \in \mathfrak{S}_r} \sum_{k=0}^r (-1)^k \begin{bmatrix}
r \\ k
\end{bmatrix} 
E_i(w_{\sigma(1)}) \cdots E_i(w_{\sigma(k)})E_j(z)E_i(w_{\sigma(k+1)}) \cdots E_i(w_{\sigma(r)}) = 0,
\end{equation}
\begin{equation}\label{loopserre2}
\sum_{\sigma \in \mathfrak{S}_r} \sum_{k=0}^r (-1)^k \begin{bmatrix}
r \\ k
\end{bmatrix} 
F_i(w_{\sigma(1)}) \cdots F_i(w_{\sigma(k)})F_j(z)F_i(w_{\sigma(k+1)}) \cdots F_i(w_{\sigma(r)}) = 0
\end{equation}
where $\delta(z)= \sum_{l \in \Z}z^l$. The equivalence of relations between generators and between currents can be easily verified by direct computation(see c.f. \cite{Her05}).

\subsection{The non-standard positive half}\label{triang}
Introduce other elements $\xi^{\pm}_{i,s}, \theta^{\pm}_{i,s}, \chi^{\pm}_{i,s}$ for $1 \leq i \leq n, s \geq 0$ via the formal functions:

$$\theta^{\pm}_i(z) = \sum_{s \geq 0} \theta^{\pm}_{i,s}z^{\pm s} = \mathrm{exp}\big( \pm(v^{-1} - v) \sum_{s \geq 0} H_{i,\pm s}z^{\pm s}\big),$$
$$\xi^{\pm}_i(z) = \sum_{s \geq 0} \xi^{\pm}_{i,s}z^{\pm s} = \mathrm{exp}\big( \sum_{s \geq 0} \frac{H_{i,\pm s}}{[s]_v}z^{\pm s}\big),$$
$$\chi^{\pm}_i(z) = \sum_{s \geq 0} \chi^{\pm}_{i,s}z^{\pm s} = \mathrm{exp}\big( -\sum_{s \geq 0} \frac{H_{i,\pm s}}{[s]_v}z^{\pm s}\big).$$

Those $H_{i,s}$'s are commutative and we can identify the subalgebras generated by $\{H_{i,s} | s \geq 0 \}$, for each $1 \leq i \leq n$, with the ring of symmetric functions $\Lambda$ by sending $H_{i,s}$ to the power sum symmetric functions $p_s$. With such identification, those elements $\xi^+_{i,s}, \chi^+_{i,s}$ and $\theta^+_{i,s}$ are nothing but(with certain normalization) the complete symmetric functions $h_s$, elementary symmetric functions $e_s$ and $q_s$ in \cite{Mac95} respectively. For any partition $\lambda$, we will denote by $\mathbf{b}_{i,\lambda}$ the element corresponding to the Schur function $s_\lambda$. We have the following property:

\begin{prop}
The sets $\{ \xi^+_{i,s} \; | \; 1 \leq i \leq n, \; s \geq 0 \}$, $\{ \chi^+_{i,s} \; | \; 1 \leq i \leq n, \; s \geq 0 \}$ or $\{ \theta^+_{i,s} \; | \; 1 \leq i \leq n, \; s \geq 0 \}$ generate the same subalgebra $U^+_v(\mathcal{L}\mathfrak{h})$ of $U_v(\mathcal{L}\mathfrak{g})$ as $\{ H_{i,+s} \; | \; 1 \leq i \leq n, \; s \geq 0 \}$. Moreover, we have the following relations:
\begin{equation}
\xi^+_i(z)\chi^+_i(z) = 1, \; \theta^+_i(z) = \xi^+_i(v^{-1}z)\chi^+_i(vz).
\end{equation}
\end{prop}

Now, let $U_v^+(\mathcal{L}\mathfrak{g})$ be the subalgebra of $U_v(\mathcal{L}\mathfrak{g})$ generated by $\{E_{i,l}, H_{i,s} \; | \; i \in I, l \in Z, s \geq 0 \}$. The triangular decomposition(\cite[Section 3.3]{Her05}) of $U_v(\mathcal{L}\mathfrak{g})$ gives rise to a decomposition $U_v^+(\mathcal{L}\mathfrak{g}) \simeq U_v(\mathcal{L}\mathfrak{n}) \otimes U^+_v(\mathcal{L}\mathfrak{h})$ via the multiplication map. If $\mathfrak{g}$ is a simple Lie algebra, then $U_v^+(\mathcal{L}\mathfrak{g})$ is nothing but the Drinfeld's positive part of the quantum affine algebra $U_v(\hat{\mathfrak{g}})$(without central charge).

\section{Bar-involution}

\subsection{Drinfeld's new comultiplication}

One can define a coproduct on $U_v(\mathcal{L}\mathfrak{g})$ (known as the \textit{Drinfeld's new coproduct}) by the following formula which define on $U_v(\mathcal{L}\mathfrak{g})$ the structure of a topological bialgebra:

$$\Delta(E_{i,n}) = E_{i,n} \otimes 1 + \sum_{t \geq 0} \phi^+_{i,t} \otimes E_{i,n-t} $$
$$\Delta(F_{i,n}) = 1 \otimes F_{i,n} + \sum_{t \geq 0} F_{i,n-t} \otimes \phi^{-}_{i,-t} $$
$$ \Delta(\phi^{\pm}_{i,\pm m}) = \sum_{0 \leq t \leq m} \phi^{\pm}_{i,\pm (m-t)} \otimes \phi^{\pm}_{i, \pm t} $$

\subsection{Admissible forms}

To introduce the bilinear form, we first define an algebra structure on $U_v^+(\mathcal{L}\mathfrak{g}) \hat{\otimes} U_v^+(\mathcal{L}\mathfrak{g})$ by 
\begin{equation*}
(x_1 \otimes x_2)(y_1 \otimes y_2) = v^{-(\mathrm{wt}x_2, \mathrm{wt}y_1)} x_1 y_1 \otimes x_2 y_2,
\end{equation*}
where $x_t, y_t$ for $t = 1,2$ are homogeneous. Let $\Delta' : U_v^+(\mathcal{L}\mathfrak{g}) \to U_v^+(\mathcal{L}\mathfrak{g}) \hat{\otimes} U_v^+ (\mathcal{L}\mathfrak{g})$ be the $\CC(v)$-algebra homomorphism defined by extending $\Delta'(E_{i,n}) = E_{i,n} \otimes 1 + \sum_{t \geq 0} \theta^+_{i,t} \otimes E_{i,n-t}$, $\Delta'(\theta^+_{i,m}) = \sum_{0 \leq t \leq m} \theta^+_{i,m-t} \otimes \theta^+_{i, t}$ for $i \in I, n \in \Z, m \in \Z_{\geq 0}$. We have a nondegenerate symmetric bilinear form $( \  , \  ): U_v^+(\mathcal{L}\mathfrak{g}) \times U_v^+(\mathcal{L}\mathfrak{g}) \to \CC $ defined by 
\begin{equation*}
( E_{i,k}, E_{j,l}) =  \frac{\delta_{i,j} \delta_{k,l}}{v^{-2}-1}, \ (H_{i,m}, H_{i,n}) = \delta_{i,j} \delta_{m,n} \frac{[2m]}{m(v^{-1} - v)}, \ (E_{i,k}, H_{j,m}) = 0. 
\end{equation*}
It is known(\cite{Gro07}) that this bilinear form is a \textit{Hopf pairing} with respect to $\Delta'$, i.e., $(x,yz) = (\Delta'(x), y \otimes z)$ where $(x\otimes y, z \otimes w) := (x ,z)(y,w)$. For any $x \in U_v^+(\mathcal{L}\mathfrak{g})$, we write $\Delta'(x) = \sum x_{(1)} \otimes x_{(2)}$. 

\vspace{.1in}

Let $U_v(\mathcal{L}\mathfrak{n})$ be the subalgebra of $U_v^+(\mathcal{L}\mathfrak{g})$ generated by $\{E_{i,l}| i \in I, l \in \Z \}$. For $i \in I$ and $ n \in \Z$, we define the linear operator $F'_{i,n}$ on $U_v(\mathcal{L}\mathfrak{n})$ by 
\begin{equation*}
F'_{i,n}(x) = \sum (v^{-2}-1)(E_{i,n}, x_{(1)})x_{(2)}
\end{equation*}
for any $x \in U_v(\mathcal{L}\mathfrak{n})$. We also let $E_{i,n}$ act on $U_v(\mathcal{L}\mathfrak{n})$ by left multiplication. By definition, the symmetric bilinear form above is $admissible$(in the sense of \cite{LuBook}, i.e. $(x,F'_{i,n}y) = (1-v_i^{-2})(E_{i,n}x,y)$) when we restrict it to $U_v(\mathcal{L}\mathfrak{n})$ and we have the following relation:
\begin{lem}\label{lemma1}
For any $i,j \in I$ and $m,l \in \Z$, 
\begin{equation}
F'_{i,l+1} F'_{j,m} - v^{b_{i,j}} F'_{j,m} F'_{i,l+1} = v^{b_{i,j}} F'_{i,l} F'_{j,m+1} - F'_{j,m+1} F'_{i,l}.
\end{equation}
For $i \neq j$, $r = 1-b_{i,j}$ and all sequences of integers $l_1,\dots,l_r$,
\begin{equation}
\sum_{\sigma \in \mathfrak{S}_r} \sum_{k =0}^r (-1)^k \begin{bmatrix}
r \\ k
\end{bmatrix}_{v_i} 
F'_{i,l_{\sigma(1)}} \cdots F'_{i,w_{\sigma(k)}}F'_{j,l'}F'_{i,w_{\sigma(k+1)}} \cdots F'_{i,w_{\sigma(r)}} = 0.
\end{equation}
\end{lem}

\subsection{Kashiwara operators}
Set

\begin{equation*}
\begin{split}
&W'_{i,k} = \sum_{ n \leq k} E_{i,n} \cdot U_v(\mathcal{L}\mathfrak{n}) \\
&Z'_{i,k} = \bigcap_{ n \leq k} \mathrm{ker} F'_{i,n}.
\end{split}
\end{equation*}
It is clear that, for each $i \in I$ and $k \in \Z$, we have a direct sum decomposition $U_v(\mathcal{L}\mathfrak{n}) \cong W'_{i,k} \oplus Z'_{i,k}$ as vector spaces. 

\begin{prop}
For any $k \in \Z$ and $i \in I$, $Z'_{i,k}$ is a subalgebra of $U_v(\mathcal{L}\mathfrak{n})$.
\end{prop}

\noindent
\textit{Proof.}
Fix $k \in \Z, i \in I$ and $n \leq k$. Let $x,y \in Z'_{i,k}$, we have 
\begin{equation*}
\begin{split}
F'_{i,n}(xy) &= \sum v^{-(\mathrm{wt}x_{(2)}, \mathrm{wt}y_{(1)})}(v^{-2}-1)(E_{i,n},x_{(1)}y_{(1)})x_{(2)}y_{(2)} \\
&=  \sum v^{-(\mathrm{wt}x_{(2)}, \mathrm{wt}y_{(1)})}(v^{-2}-1)\{ (E_{i,n},x_{(1)})(1,y_{(1)})x_{(2)}y_{(2)} \\ 
& \; + (1,x_{(1)})(E_{i,n},y_{(1)})x_{(2)}y_{(2)} + \sum_{t > 0} (\theta^+_{i,t},x_{(1)})(E_{i,n-t},y_{(1)})x_{(2)}y_{(2)} \} \\
&= \sum (v^{-2}-1)(E_{i,n},x_{(1)})x_{(2)}y + \sum (v^{-2}-1)(E_{i,n}, y_{(1)})xy_{(2)} \\
& \; +  \sum v^{-(\mathrm{wt}x_{(2)}, \mathrm{wt}y_{(1)})}(v^{-2}-1)\sum_{t > 0} (\theta^+_{i,t},x_{(1)})x_{(2)}(E_{i,n-t},y_{(1)})y_{(2)} \\
&= F'_{i,n}(x)y + xF'_{i,n}(y) +  \sum v^{-(\mathrm{wt}x_{(2)}, \mathrm{wt}y_{(1)})}(v^{-2}-1)\cdot \\ 
&\; \cdot \sum_{t > 0} (\theta^+_{i,t},x_{(1)})x_{(2)}(E_{i,n-t},y_{(1)})y_{(2)} \\
&=  \sum v^{-(\mathrm{wt}x_{(2)}, \mathrm{wt}y_{(1)})}(v^{-2}-1)\sum_{t > 0} (\theta^+_{i,t},x_{(1)})x_{(2)}(E_{i,n-t},y_{(1)})y_{(2)},
\end{split}
\end{equation*}
where for the second equality we use the fact that the bilinear form is a Hopf pairing and the last one holds since $x,y \in Z'_{i,k}$. To show the remaining term vanish, let us fix $x_{(1)},x_{(2)}$. Note that $(E_{i,n-t},y_{(1)}) \neq 0$ only if $ y_{(1)}$ is a scalar multiple of $E_{i,l}$ for some $l \in \Z$. Hence all the $(\mathrm{wt}x_{(2)}, \mathrm{wt}y_{(1)})$ are equal, say $C(x_{(2)})$, and we have 
\begin{equation*}
\begin{split}
 & \sum v^{-(\mathrm{wt}x_{(2)}, \mathrm{wt}y_{(1)})}(v^{-2}-1)\sum_{t > 0} (\theta^+_{i,t},x_{(1)})x_{(2)}(E_{i,n-t},y_{(1)})y_{(2)} \\
=& \sum C(x_{(2)}) (v^{-2}-1))\sum_{t > 0} (\theta^+_{i,t},x_{(1)})x_{(2)}\{\sum (E_{i,n-t},y_{(1)})y_{(2)} \} \\
=& \sum C(x_{(2)}) \sum_{t > 0} (\theta^+_{i,t},x_{(1)})x_{(2)} F'_{i,n-t}(y) = 0
\end{split}
\end{equation*}

\qed

\begin{lem}\label{commutative}
For each $i,j \in I$ and $m,n \in \Z$, we have the following relations
\begin{equation*}
F'_{i,m}E_{j,n} - v^{-b_{i,j}}E_{j,n}F'_{i,m} = \delta_{i,j}\delta_{n,m} + \sum_{t > 0} v^{-tb_{i,j}}(v^{-b_{i,j}} - v^{b_{i,j}})E_{j,n-t}F'_{i,m-t}.
\end{equation*}
\end{lem}

\noindent
\textit{Proof.} 
By direct computation we obtain the formula
\begin{equation}
\theta^+_{i,l} E_{j,n} = E_{j,n}\theta^+_{i,l} + \sum_{t = 1}^{l} v^{-(t-1)b_{i,j}}(v^{-b_{i,j}} - v^{b_{i,j}})E_{j,n+t}\theta^+_{i,l-t}.
\end{equation}
For $x \in U_v(\mathcal{L}\mathfrak{n})$, we calculate 
\begin{equation*}
\begin{split}
F'_{i,m}E_{j,n}(x) &= F'_{i,m}(E_{j,n}x) \\
&= (v^{-2}-1)(E_{i,m}, E_{j,n})x \\
&+ (v^{-2}-1)\sum v^{-(\mathrm{wt}E_{j,n}, \mathrm{wt}x_{(1)} )}(E_{i,m}, \theta^+_{j,0}x_{(1)})E_{j,n}x_{(2)} \\
&+ (v^{-2}-1)\sum \sum_{t >0} v^{-(\mathrm{wt}E_{j,n-t}, \mathrm{wt}x_{(1)} )}(E_{i,m}, \theta^+_{j,t}x_{(1)})E_{j,n-t} x_{(2)}
\end{split}
\end{equation*}
All the $(E_{i,m}, \theta^+_{j,t}x_{(1)})$ must vanish but $x_{(1)}$ is a scalar multiple of $E_{i,l}$ for some $l \in \Z$. Hence we have 
\begin{equation*}
\begin{split}
F'_{i,m}E_{j,n}(x) &= (v^{-2}-1)(E_{i,m}, E_{j,n})x \\
& + v^{-b_{i,j}}(v^{-2}-1) \sum (E_{i,m}, x_{(1)}) E_{j,n} x_{(2)} \\
& + \sum_{t > 0} (v^{-2}-1) \sum v^{-t b_{i,j}}(v^{-b_{i,j}} - v^{b_{i,j}}) (E_{i,m}, x_{(1)+t}) E_{j, n-t} x_{(2)} \\
&= \delta_{i,j} \delta_{m,n} + v^{-b_{i,j}} E_{j,n}F'_{i,m}(x) \\
& + \sum_{t > 0} v^{-t b_{i,j}}(v^{-b_{i,j}} - v^{b_{i,j}}) (E_{i,m-t}, x_{(1)}) E_{j,n-t} x_{(2)} \\
&=  \delta_{i,j} \delta_{m,n} + v^{-b_{i,j}} E_{j,n}F'_{i,m}(x) \\
& + \sum_{t > 0} v^{-tb_{i,j}}(v^{-b_{i,j}} - v^{b_{i,j}})E_{j,n-t}F'_{i,m-t}(x),
\end{split}
\end{equation*}
where $x_{(1)+t} := cE_{i,l+t}$ if $x_{(1)} = c E_{i,l}$ for some scalar $c \in \mathbb{C}$.

\qed

\vspace{.1in}

Fix an $i \in I$. For any $n \in \Z$, $k \leq n-1$ and $x \in Z'_{i,n-1}$, by Lemma \ref{commutative}, we have
\begin{equation*}
F'_{i,k}(E_{i,n} \cdot x) = v^2E_{i,n}F'_{i,k}(x) + \sum_{t>0} v^{-2t}(v^{-2} - v^2) E_{i,n-t} F'_{i,n-t}(x) = 0.
\end{equation*}
Hence the space $Z'_{i,n-1}$ is stable under the action of $E_{i,n}$. Moreover, by Lemma \ref{lemma1}, $Z'_{i,n-1}$ is also stable under the action of $F'_{i,n}$. Let us fix also an $n \in \Z$ now. It follows from the Lemma \ref{commutative} again that if we restrict the operators $E_{i,n}, F'_{i,n}$ to $Z'_{i,n-1}$, we have the so-called $q$-Boson relation:
\begin{equation}
F'_{i,n} E_{i,n} - v^{-2} E_{i,n} F_{i,n} = 1.
\end{equation} 
For $s \in \Z$, let $E^{(s)}_{i,n} : Z'_{i,n-1} \to Z'_{i,n-1}$ be defined as $\frac{E^s_{i,n}}{[s]_v!}$ if $s \geq 0$ and as $0$ if $s < 0$. We can deduce by induction on $s$ from the $q$-Boson relation that:
\begin{equation}
F'_{i,n} E^{(s)}_{i,n} = v^{-2s}E^{(s)}_{i,n}F'_{i,n} + v^{-(s-1)}E_{i,n}^{(s-1)}.
\end{equation}
For any $t \geq 0$, consider the operator
\begin{equation}
\Pi_{i,n,t}  = \sum_{s \geq 0} (-1)^sv^{-s(s-1)/2} E^{(s)}_{i,n} F'^{s+t}_{i,n} : Z'_{i,n-1} \to Z'_{i,n-1}.
\end{equation}
Clearly $F'_{i,n}$ acts on $Z'_{i,n-1} $ locally nilpotently, $\Pi_{i,n,t}$ is well-defined and hence $Z'_{i,n-1}$ is an object of the category $\mathcal{D}_i$ defined in \cite[Chapter 16]{LuBook}. We have the following properties
\begin{prop}\cite[Lemma 16.1.2]{LuBook}
\begin{enumerate}
\item[(1)] We have $F'_{i,n}\Pi_{i,n,t} = 0$ for all $t \geq 0$.
\item[(2)] We have $\sum_{t \geq 0} v^{t(t-1)/2} E_{i,n}^{(t)} \Pi_{i,n,t} =1$.
\item[(3)] We have a direct sum decomposition $Z'_{i,n-1} = \bigoplus_{N \geq 0} Z'_{i,n-1}(N) $ as a vector space, where for $N \geq 0$, $Z'_{i,n-1}(0) = \{ x \in Z'_{i,n-1} \mid F'_{i,n}(x) = 0 \}$ and $Z'_{i,n-1}(N) = E^{(N)}_{i,n} Z'_{i,n-1}(0)$. Moreover, the map $E^{(N)}_{i,n}$ restricts to an isomorphism of vector spaces $Z'_{i,n-1}(0) \simeq Z'_{i,n-1}(N)$.
\item[(4)] $E_{i,n}: Z'_{i,n-1} \to Z'_{i,n-1}$ is injective. 
\end{enumerate}
\end{prop}
Thus we can define the linear maps $\widetilde{E}_{i,n},\widetilde{F}_{i,n} : Z'_{i,n-1} \to Z'_{i,n-1}$ by 
\begin{equation}
\begin{split}
\widetilde{E}_{i,n}(E_{i,n}^{(N)}y) & = E^{(N+1)}_{i,n} y , \\
\widetilde{F}_{i,n}(E_{i,n}^{(N)}y) & = E^{(N-1)}_{i,n} y
\end{split}
\end{equation}
for any $y \in Z'_{i,n-1}(0)$ and extend them to $U_v(\mathcal{L}\mathfrak{n}) = W'_{i,n-1} \oplus Z'_{i,n-1}$ by sending $x \in W'_{i,n-1}$ to $0$. These operators $\widetilde{E}_{i,n},\widetilde{F}_{i,n}$ can be thought as the loop's analogue of \textit{Kashiwara's operators}. We might also extend the operators $\widetilde{E}_{i,n},\widetilde{F}_{i,n}$ to $U^+_v(\mathcal{L}\mathfrak{g})$ by defining $\widetilde{E}_{i,n}(x) = \widetilde{E}_{i,n}(y)h$ for any $x = yh \in U^+_v(\mathcal{L}\mathfrak{g}) = U_v(\mathcal{L}\mathfrak{n})U^+_v(\mathcal{L}\mathfrak{h})$ with respect to the triangular decomposition in Section \ref{triang}. 

\vspace{.1in}

\noindent
\textit{Remark}. The algebra $B_v(\mathcal{L}\mathfrak{g})$ generated by $E_{i,n}, F'_{i,n}$ for all $i \in I$ and $n \in \Z$ can be thought as the \textit{affinization} of the $q$-Boson algebra $B_v(\mathfrak{g})$.

\vspace{.1in}

Let $\mathcal{A}$ be the subring of $\CC(v)$ consisting of rational functions without pole at $v = 0$. Let $\mathcal{L}'$ be the $\mathcal{A}$-submodule of $U_v(\mathcal{L}\mathfrak{n})$ generated by $\widetilde{E}_{i_1,n_1} \dots \widetilde{E}_{i_l,n_l}\cdot 1$ with $l \in \mathbb{N}, i_1,\dots,i_l \in I, n_1,\dots,n_l \in \mathbb{Z}^l$ and let $\mathcal{B}'$ be the set of their images in $\mathcal{L}' / v \mathcal{L}'$. Similarly, we define $\mathcal{L}$ be the $\mathcal{A}$-submodule of $U^+_v(\mathcal{L}\mathfrak{g})$ generated by $\widetilde{E}_{i_1,n_1} \dots \widetilde{E}_{i_l,n_l}\cdot \mathbf{b}_{j,\lambda}$ with $j \in I$ and all partitions $\lambda$ and let $\mathcal{B}$ be the set of their images in $\mathcal{L} / v \mathcal{L}$.

\begin{conj}\label{mainconj}
We have
\begin{enumerate}
\item $\mathcal{L}$ is a free $\mathcal{A}$-module such that $\CC(v) \otimes_{\mathcal{A}} \mathcal{L} \cong U_v^+(\mathcal{L}\mathfrak{g})$.
\item $\mathcal{B}$ is a basis for $\mathcal{L} / v \mathcal{L}$.
\item $\mathcal{L}$ is stable by $\widetilde{E}_{i,l},\widetilde{F}_{i,l}$ for all $i \in I$ and $l \in \mathbb{Z}$.
\item $\widetilde{E}_{i,l}\mathcal{B}, \widetilde{F}_{i,l}\mathcal{B} \subset \mathcal{B} \cup \{ 0 \}$ for all $i \in I$ and $l \in \mathbb{Z}$.
\item For $b,b' \in \mathcal{B}$ one has $b' = \widetilde{E}_{i,l}b \Longleftrightarrow \widetilde{F}_{i,l}b' = b$. 
\end{enumerate}
\end{conj}

\subsection{The completion}
In this section we are going to define a completion(similar to the Harder-Harasimhan completion) $\widehat{U}^+_v(\mathcal{L}\mathfrak{g})$ of $U_v^+(\mathcal{L}\mathfrak{g})$. 

\vspace{.1in}

Let us write $U_v^+(\mathcal{L}\mathfrak{g}) = \bigoplus_{\alpha \in \widehat{Q}^+}U_v^+(\mathcal{L}\mathfrak{g})[\alpha]$ as the root spaces decomposition. We define a \textit{slope function} $\mu: \widehat{Q}^+ \to \mathbb{Q} \cup \mathbb{Z}$ by $\mu(\alpha) = \frac{d}{|\alpha_0|}$ for any $\alpha = (\alpha_0,d) \in  Q^+ \times \mathbb{Z} = \widehat{Q}^+$. Fix a root $\alpha \in \widehat{Q}^+$ and for any $m \in \mathbb{Z}$ we set
\begin{equation*}
W_m[\alpha] = \sum_{\beta \leq \alpha, \mu(\beta) \leq m} U_v^+(\mathcal{L}\mathfrak{g})[\beta] \cdot U_v^+(\mathcal{L}\mathfrak{g})[\alpha - \beta] \subset U_v^+(\mathcal{L}\mathfrak{g})[\alpha]
\end{equation*}
so that we have $U_v^+(\mathcal{L}\mathfrak{g})[\alpha] = W_m[\alpha] \oplus Z_m[\alpha]$ where $Z_m[\alpha] := U_v^+(\mathcal{L}\mathfrak{g})[\alpha]/ W_m[\alpha]$. For any pair $m \geq n$ the canonical embedding $W_n[\alpha] \to W_m[\alpha]$ induces a commutative diagram
\begin{equation*}
\xymatrix
{
U_v^+(\mathcal{L}\mathfrak{g})[\alpha] / W_n[\alpha] \ar[r]^-{\pi_n} \ar[d] & Z_n[\alpha] \ar[d]^-{\phi_{m,n}} \\
U_v^+(\mathcal{L}\mathfrak{g})[\alpha] / W_m[\alpha]  \ar[r]^-{\pi_m} & Z_m[\alpha]  
}
\end{equation*}
Obviously $(Z_m[\alpha], \phi_{m,n})$ forms a projective system and we can define 
\begin{equation*}
\widehat{U}^+_v(\mathcal{L}\mathfrak{g})[\alpha] := \varprojlim Z_{m}[\alpha].
\end{equation*}
Each $Z_m[\alpha]$ is finite dimensional and since $U^+_v(\mathcal{L}\mathfrak{g})[\alpha] = \cup_m Z_m[\alpha]$ we may view $\widehat{U}^+_v(\mathcal{L}\mathfrak{g})[\alpha]$ as the set of infinite sums $\sum c_x x$ with $c_x \in \CC(v), \mathrm{wt}(x) = \alpha$. For the sake of convenience we also denote by $\textit{jet}_m$ the canonical morphism $\widehat{U}^+_v(\mathcal{L}\mathfrak{g})[\alpha] \to Z_m[\alpha]$. By the universal property of the projective limit there is an injective linear map $U^+_v(\mathcal{L}\mathfrak{g})[\alpha] \to \widehat{U}^+_v(\mathcal{L}\mathfrak{g})[\alpha]$, and since the map $U^+_v(\mathcal{L}\mathfrak{g})[\alpha] \to W_m[\alpha]$ splits, we may consider $W_m[\alpha]$ as a subspace of $\widehat{U}^+_v(\mathcal{L}\mathfrak{g})[\alpha]$ via the inclusion $W_m[\alpha] \to U^+_v(\mathcal{L}\mathfrak{g})[\alpha] \to \widehat{U}^+_v(\mathcal{L}\mathfrak{g})[\alpha]$. So the projection $\textit{jet}_m : \widehat{U}^+_v(\mathcal{L}\mathfrak{g})[\alpha] \to Z_m[\alpha]$ is an idempotent morphism. Let us denote $r_m = 1 - \textit{jet}_m$. Then any element $u \in \widehat{U}^+_v(\mathcal{L}\mathfrak{g})[\alpha]$ can uniquely be written as $\textit{jet}_m(u) + r_m(u)$, where $\textit{jet}_m(u) \in Z_m[\alpha]$ and $\textit{jet}_m(r_m(u)) = 0$. Using this formalism, the space $U^+_v(\mathcal{L}\mathfrak{g})[\alpha]$ as a subset of $\widehat{U}^+_v(\mathcal{L}\mathfrak{g})[\alpha]$ can be identified with the set of those sequence $u = (u_m)$ for which $r_m(u_m) = 0$ for $m << 0$. Finally, we define 
\begin{equation}
\widehat{U}^+_v(\mathcal{L}\mathfrak{g}) = \bigoplus_{\alpha \in \widehat{Q}^+} \widehat{U}^+_v(\mathcal{L}\mathfrak{g})[\alpha].
\end{equation}

\begin{prop}
$\widehat{U}^+_v(\mathcal{L}\mathfrak{g})$ is an associative algebra.
\end{prop}

\noindent
\textit{Proof.}
Let us show that, for any $\alpha, \beta \in \widehat{Q}^+$, the multiplication 
\begin{equation*}
\widehat{U}^+_v(\mathcal{L}\mathfrak{g})[\alpha] \otimes \widehat{U}^+_v(\mathcal{L}\mathfrak{g})[\beta] \to \widehat{U}^+_v(\mathcal{L}\mathfrak{g})[\alpha + \beta]
\end{equation*}
defined by component-wise multiplication is well-defined. We need to show that, for each $n \in \Z $, $\alpha,\beta \in \widehat{Q}^+$ and $a \in \widehat{U}^+_v(\mathcal{L}\mathfrak{g})[\alpha], b \in \widehat{U}^+_v(\mathcal{L}\mathfrak{g})[\beta]$, there exist $l, l' \ll 0 $ such that $ab \equiv a_{m} b_{m'}$(mod $W_n[\alpha + \beta]$) for all $m < l, m' < l'$. In other words, we need to show that there exist $l,l' \ll 0$ such that $W_m[\alpha] b \subset W_n[\alpha+\beta]$ and $a W_{m'}[\beta] \subset W_n[\alpha+\beta]$ for all $a \in U^+_v(\mathcal{L}\mathfrak{g})[\alpha], b \in U^+_v(\mathcal{L}\mathfrak{g})[\beta]$ and $m < l, m' < l'$. Clearly the first inclusion follows from the definition of $W_m[\alpha]$. To show the later one, let us set
\begin{equation*}
k(a) = \min \{ k \in \mathbb{Z} \mid a \in W_k[\alpha] \}. 
\end{equation*}
If $k(a) \leq n$, then $aW_{m'}[\beta] \subset W_n[\alpha+\beta]$ for all $m' \in \Z$. Suppose now $k(a) > n$. Using the quadratic relation (\ref{relation1}), we have $aW_{m'}[\beta] \subset W_{m''}[\alpha + \beta]$ where $m'' = m' + |\alpha_0|$ and any $l' \leq n-|\alpha_0|$ will be suitable. Thus the multiplication is well-defined. The associativity follows from the associativity of the multiplication of $U^+_v(\mathcal{L}\mathfrak{g})$.
\qed

\vspace{.2in}

\subsection{Bar-involution}
Let us define the semi-linear map $ \varphi : U_v^+(\mathcal{L}\mathfrak{g}) \to \widehat{U}_v^+(\mathcal{L}\mathfrak{g})$ on the generators by
\begin{equation*}
\begin{split}
&v \mapsto \varphi( v)= v^{-1} \\
&\xi^+_{i,s} \mapsto \varphi(\xi^+_{i,s}) = \xi^+_{i,s} \\
&E_{i,l} \mapsto \varphi(E_{i,l}) = E_{i,l} + \sum_{t=1}^\infty (-1)^t \sum_{l_1, \cdots ,l_t >0 }(v^{-l_t} - v^{l_t})v^{-\sum_{j=1}^{t-1}l_j} E_{i,l-\sum_{j=1}^t l_j} \prod_{j=1}^t \xi^+_{i,l_j}
\end{split}
\end{equation*}
for any $i \in I ,l \in \Z, s \geq 0$ and extend it to $U^+_v(\mathcal{L}\mathfrak{g})$ by setting $\varphi(xy)= \varphi(x) \varphi(y)$ for any $x,y \in U^+_v(\mathcal{L}\mathfrak{g})$. Clearly the map $\varphi$ preserves the weight spaces and its image lies in the completion $\widehat{U}^+_v(\mathcal{L}\mathfrak{g})$. One can easily check that $\varphi(\varphi(E_{i,l})) = E_{i,l}$ for any $i \in I$ and $l \in \Z$. This section is devoted to prove the following:

\begin{theo}
The involution $\varphi: U^+_v(\mathcal{L}\mathfrak{g}) \to \widehat{U}^+_v(\mathcal{L}\mathfrak{g})$ is a ring involution and extends to an involution $\varphi : \widehat{U}^+_v(\mathcal{L}\mathfrak{g}) \to \widehat{U}^+_v(\mathcal{L}\mathfrak{g})$. 
\end{theo} 

\noindent
\textit{Proof.} We will use the generating series to verify the defining relations. Set
\begin{equation}
\begin{split}
&E^0_i(z) = \xi^+_i(\frac{z}{v}) \\
&E^1_i(z) = E_i(z)E^0_i(v^2 z)
\end{split}
\end{equation}
Then $\varphi(E^0_i(z)) = E^0_i(v^2 z)$ and $\varphi(E^1_i(z)) = E^1_i(z)$. From the definition of $E^1_i(z)$ and its bar invariance, we get
\begin{equation*}
\begin{split}
\varphi(E_i(z)) &= E^1_i(z) E^0_i(v^2 z)^{-1} = E^1_i(z) \chi^+_i(v^{-1}z) \\ &= E_i(z) \xi^+_i(vz) \chi^+_i(v^{-1}z) \\
&= E_i(z) (\theta^+_i(z))^{-1}.
\end{split}
\end{equation*}
We first compute $\varphi(E_i(z))\varphi(E_j(w))$:

\begin{equation}\label{propeq}
\begin{split}
&\varphi(E_i(z))\varphi(E_j(w)) \\ 
&= E_i(z)\theta^+_i(z)^{-1}E_j(w)\theta^+_j(w)^{-1} \\
&= \frac{1-\frac{z}{w}v^{-b_{i,j}}}{1-\frac{z}{w}v^{b_{i,j}}}E_i(z)E_j(w)\theta^+_i(z)^{-1}\theta^+_j(w)^{-1} \\
&= v^{-b_{i,j}} E_j(w)E_i(z)\theta^+_i(z)^{-1}\theta^+_j(w)^{-1}  .
\end{split}
\end{equation}

To check the quadratic relation (\ref{loop2}), 
\begin{equation*}
\begin{split}
& \varphi \big( (wv^{b_{i,j}}-z)E_i(z)E_j(w) \big) \\
= & \varphi((wv^{b_{i,j}}-z)) \varphi(E_i(z)) \varphi(E_j(w)) \\
= &(wv^{-b_{i,j}}-z)  \frac{1-\frac{z}{w}v^{-b_{i,j}}}{1-\frac{z}{w}v^{b_{i,j}}}E_i(z)E_j(w)\theta^+_i(z)^{-1}\theta^+_j(w)^{-1} \\
= &v^{-2b_{i,j}}(w v^{b_{i,j}}-z) E_i(z)E_j(w)\theta^+_i(z)^{-1}\theta^+_j(w)^{-1} \\
= &v^{-2b_{i,j}}(w-zv^{b_{i,j}})E_j(w)E_i(z)\theta^+_i(z)^{-1}\theta^+_j(w)^{-1} \\
= &v^{-2b_{i,j}}(w-zv^{b_{i,j}})\frac{1-\frac{w}{z}v^{b_{i,j}}}{1-\frac{w}{z}v^{-b_{i,j}}} E_j(w)\theta^+_j(w)^{-1}E_i(z)\theta^+_i(z)^{-1} \\
= &(w-z v^{-b_{i,j}}) E_j(w)\theta^+_j(w)^{-1}E_i(z)\theta^+_i(z)^{-1} \\
= &\varphi((w-z v^{b_{i,j}})) \varphi(E_j(w)) \varphi(E_i(z)) \\
= &\varphi \big( (w-z v^{b_{i,j}})E_j(w)E_i(z) \big) .
\end{split}
\end{equation*}
\vspace{.1in}

For the Serre relation (\ref{loopserre1}), by using (\ref{propeq}) we have
\begin{equation*}
\begin{split}
& \varphi(E_i(w_{\sigma(1)}) )\cdots \varphi(E_i(w_{\sigma(k)}) ) \varphi(E_j(z) ) \varphi(E_i(w_{\sigma(k+1)}) )\cdots \varphi(E_i(w_{\sigma(r)}) ) \\
&=  E_i(w_{\sigma(1)})\theta^+_i(w_{\sigma(1)})^{-1} \cdots E_i(w_{\sigma(k)})\theta^+_i(w_{\sigma(k)})^{-1} E_j(z)\theta^+_j(z)^{-1} \\ & \cdot E_i(w_{\sigma(k+1)})\theta^+_i(w_{\sigma(k+1)})^{-1} \cdots E_i(w_{\sigma(r)})\theta^+_i(w_{\sigma(r)})^{-1} \\
&= v^{-r(r-1)-rb_{i,j}} E_i(w_{\sigma(r)}) \cdots E_i(w_{\sigma(k+1)})E_j(z)E_i(w_{\sigma(k)})\cdots E_i(w_{\sigma(1)}) \cdot \\
& \cdot \prod_{l =1}^r \theta^+_i(w_l)^{-1} \theta^+_j(z)^{-1}
\end{split}
\end{equation*}
Hence
\begin{equation*}
\begin{split}
& \sum_{\sigma \in \mathfrak{S}_r} \sum_{k=0}^r (-1)^k \begin{bmatrix}
r \\ k
\end{bmatrix} \cdot \\
& \cdot \varphi(E_i(w_{\sigma(1)}) ) \cdots \varphi(E_i(w_{\sigma(k)})) \varphi(E_j(z)) \varphi(E_i(w_{\sigma(k+1)})) \cdots \varphi(E_i(w_{\sigma(r)})) \\
&= v^{-r(r-1)-b_{i,j}}\big( \sum_{\sigma \in \mathfrak{S}_r} \sum_{k=0}^r (-1)^k \begin{bmatrix}
r \\ k
\end{bmatrix} 
E_i(w_{\sigma(r)}) \cdots E_i(w_{\sigma(k+1)})E_j(z)E_i(w_{\sigma(k)})\\
& \cdots E_i(w_{\sigma(1)}) \big) \cdot \prod_{l =1}^r \theta^+_i(w_l)^{-1} \theta^+_j(z)^{-1} \\
&= 0
\end{split}
\end{equation*}
Thus $\varphi : U^+_v(\mathcal{L}\mathfrak{g}) \to \widehat{U}^+_v(\mathcal{L}\mathfrak{g})$ is indeed a ring involution. To show that it extends to $\widehat{U}^+_v(\mathcal{L}\mathfrak{g}) \to \widehat{U}^+_v(\mathcal{L}\mathfrak{g})$, we have to show the continuity of the involution, i.e., for any $n \in \Z$ and $\alpha \in \widehat{Q}^+$, we have $\varphi(W_n[\alpha]) \subset \overline{W_n[\alpha]}$, where $\overline{W_n[\alpha]}$ denotes the closure of $W_n[\alpha]$ in $\widehat{U}^+_v(\mathcal{L}\mathfrak{g})$. However, the involution $\varphi$ preserves the weights and slopes, the continuity follows from the density of $W_n[\alpha]$ in $\overline{W_n[\alpha]}$. 

\qed

\subsection{Examples and conjectures}
t is clear that $W'_{i,k}[\alpha] = W'_{i,k} \cap U_v(\mathcal{L}\mathfrak{n})[\alpha] \subset W_k[\alpha]$ for all $\alpha \in \widehat{Q}^+$, $i \in I$ and all $k \in \Z$. The Kashiwara operators $\widetilde{E}_{i,k},\widetilde{F}_{i,k}$ are well-defined on $\widehat{U}_v(\mathcal{L}\mathfrak{n})$, the closure of $U_v(\mathcal{L}(\mathfrak{n}))$ in $\widehat{U}^+_v(\mathcal{L}\mathfrak{g})$ as well as on $\widehat{U}^+_v(\mathcal{L}\mathfrak{g})$ too. We set $\widehat{\mathcal{L}} : = \overline{\mathcal{L}}$ and $\widehat{\mathcal{L}}^- := \overline{\varphi(\mathcal{L})}$. Let us assume that the Conjecture \ref{mainconj} is true. Then the density of $U^+_v(\mathcal{L}\mathfrak{g})$ in $\widehat{U}^+_v(\mathcal{L}\mathfrak{g})$ implies that $\mathbb{C}(v) \otimes_{\mathcal{A}} \widehat{\mathcal{L}} \simeq \widehat{U}^+_v(\mathcal{L}\mathfrak{g})$. Moreover, the subset $\widehat{\mathcal{B}}$ of the images of lattice $\widehat{\mathcal{L}}$ in $\widehat{\mathcal{L}} / v \widehat{\mathcal{L}}$ is a basis. Thus we have the following conjecture of the existence of the \textit{canonical basis} on $\widehat{U}^+_v(\mathcal{L}\mathfrak{g})$:

\begin{conj}\label{conjecture2}
$(\widehat{U}^+_v(\mathcal{L}\mathfrak{g}), \widehat{\mathcal{L}}, \widehat{\mathcal{L}}^-)$ is a balanced triple.
\end{conj}

\vspace{.1in}

\noindent
\textbf{Example}. If $\mathfrak{g} = \mathfrak{sl}_2$, $U^+_v(\mathcal{L}\mathfrak{sl}_2)$ can be realized as the Hall algebra $\mathbf{H}_{\Coh(\mathbb{P}^1)}$ associated to the category of coherent sheaves on $\mathbb{P}^1(\mathbb{F}_l)$, where $l^{-1/2} = v$, and the completion $\widehat{U}^+_v(\mathcal{L}\mathfrak{sl}_2)$ coincides with the Harder-Narasimhan completion $\widehat{\mathbf{H}}_{\Coh(\mathbb{P}^1)}$ of $\mathbf{H}_{\Coh(\mathbb{P}^1)}$ via the HN filtration(cf. \cite[Section 2]{SchHN}). The lattice 
\[
\widehat{\mathcal{L}} = \{ E_{n_1}^{(s_1)} \cdots E_{n_l}^{(s_l)} \mathbf{b}_{\lambda} \mid n_1 > \cdots > n_l, s_1,\dots,s_l \in \mathbb{N}, \lambda \textit{ partitions}  \}
\]
is indeed a free $\mathcal{A}$-module by using the quadratic relations (\ref{relation1}) and the resulting basis $\widehat{\mathcal{B}}$ coincides with the canonical basis obtained by Schiffmann in \cite{Sch06}. Therefore the Conjecture \ref{mainconj} is true in this case.

\vspace{.2in}
\small{

\end{document}